\documentclass[article,12pt]{article}
\usepackage{setspace}
\usepackage{authblk}
\usepackage{amssymb}
\usepackage{amsmath}
\usepackage{latexsym,amssymb,amsbsy,rotating,amsxtra}
\usepackage[T1]{fontenc}
\usepackage[latin1]{inputenc}
\usepackage[english]{babel}
\usepackage{graphicx,pict2e}
\usepackage{algorithm}
\floatname{algorithm}{Algorithm}
\usepackage{natbib}

\usepackage{tikz}
\usepackage{lscape}

\newtheorem{prop}{Proposition}
\newtheorem{theo}{Theorem}

\begin{document}
\date{}

\title{A comparison of pivotal sampling and unequal probability sampling with replacement}

\author[1]{Guillaume Chauvet}
\affil[1]{ENSAI/IRMAR, Campus de Ker Lann, 35170 Bruz, France}
\author[2]{Anne Ruiz-Gazen}
\affil[2]{Toulouse School of Economics, 21 all\'ee de Brienne, 31000 Toulouse, France}

\maketitle

%\begin{center}
%\author{Guillaume Chauvet$^{(1)}$ and Anne Ruiz-Gazen$^{(2)}$ \\
%$^{(1)}$ ENSAI-IRMAR, Campus de Ker Lann, 35170 Bruz, France \\
%$^{(2)}$ Toulouse School of Economics, 21 all\'ee de Brienne, 31000 Toulouse, France}
%\end{center}
%\maketitle

\begin{abstract}
\noindent We prove that any implementation of pivotal sampling is more efficient than multinomial sampling. This yields the weak consistency of the Horvitz-Thompson estimator and the existence of a conservative variance estimator. A small simulation study supports our findings.
\end{abstract}

%\begin{keyword}
%mean-square consistency \sep sampling algorithm \sep variance estimator
%\end{keyword}

%\end{frontmatter}

\section{Introduction}

\noindent Many procedures exist for without-replacement unequal probability sampling. Pivotal sampling \citep{dev:til:98,til:11,cha:12} is a very simple sequential procedure. It satisfies strong properties of negative dependence, as proved in \citet{dub:jon:ran:07,bro:cut:rad:11,bra:jon:11,jon:12}. In particular, pivotal sampling avoids the selection of neighbouring units. This makes it particularly useful for spatial sampling, where it enables selecting samples well spread over space. A vast literature has recently focused on such applications for pivotal sampling, see for example \citet{gra:lun:sch:12,gra:rin:13,gra:saa:ene:13,gra:til:13,dic:ben:giu:esp:14,ben:pie:pos:15,dic:til:15,fat:cor:chi:pag:15,val:fer:riv:til:15}. Pivotal sampling has also found uses for longitudinal surveys \citep{ned:qua:til:09}.

\noindent For a sampling design, basic properties for estimation are that: a) the Horvitz-Thompson (HT) estimator is weakly consistent for the true total; b) the HT-estimator satisfies a central-limit theorem; c) a weakly consistent variance estimator is available for the HT-estimator. Two of these properties are tackled in this paper. We prove that any implementation of pivotal sampling is more efficient than multinomial sampling, which gives a) under a simple moment condition; see \citet[Section 5.4]{til:11} for a detailed description of multinomial sampling. It is not possible to prove c), since some second-order inclusion probabilities are zero for pivotal sampling leading to biased variance estimators. However, we prove that the Hansen-Hurvitz (HH) variance estimator \citep[equation 5.5]{til:11} provides an upper bound for the true variance, which enables to give conservative confidence intervals.

\noindent Using sufficient conditions listed in \citet{gab:90}, it has been proved that the Sampford design and the conditional Poisson sampling design \citep{gab:81,gab:84,qua:08} are more efficient than multinomial sampling. Some simple counterexamples prove that none of these sufficient conditions hold for pivotal sampling. Different tools are thus used in the present paper to obtain this property for pivotal sampling. Some basic notations are given in Section \ref{sec:not}. Ordered pivotal sampling and multinomial sampling are briefly presented in Section \ref{sec:ops}. Our main results are given in Section \ref{sec:comp}. A small simulation study supporting our theoretical results is presented in Section \ref{sec:simu}.

\section{Notations} \label{sec:not}

\noindent Consider a finite population $U$ consisting of $N$ sampling units represented by integers $k=1,\ldots,N$. Denote by $\pi=\left(\pi_1,\ldots,\pi_N\right)^{\top}$ a vector of probabilities, with $0 < \pi_k \leq 1$ for any unit $k$ in $U$, and $n=\sum_{k \in U} \pi_k$ the expected sample size. Let $p(\cdot)$ denote a sampling algorithm with parameter $\pi$, that is, such that the expected number of draws for unit $k$ in the sample equals $\pi_k$. We note $E(\cdot)$ and $V(\cdot)$ for the expectation and variance. For any variable of interest $y$, the total $t_y=\sum_{k \in U} y_k$ is unbiasedly estimated by $\hat{t}_{y}^p = \sum_{k \in S^p} \check{y}_k$ where $\check{y}_k=\pi_k^{-1} y_k$, with $S^p$ a sample selected by means of the sampling algorithm $p(\cdot)$. In case of with-replacement sampling, a unit $k$ may appear several times in $S^p$ and $\hat{t}_{y}^p$ is the \cite{han:hur:43} estimator. In case of without-replacement sampling, a unit $k$ may appear only once in $S^p$ and $\hat{t}_{y}^p$ is the \cite{hor:tho:52} estimator.

\noindent We define the cumulated inclusion probabilities for unit $k$ as $C_k=\sum_{l=1}^k \pi_l$, with $C_0=0$. The unit $k$ is cross-border if $C_{k-1} < i$ and $C_{k} \ge i$ for some integer $i=1,\ldots,n-1$. These cross-border units are denoted as $k_i$, and we note $a_i=i-C_{k_i-1}$ and $b_i=C_{k_i}-i$, for $i=1,\ldots,n-1$. We now describe a clustering of the units in $U$ which leads to a population denoted as $U_c$, so as to simplify demonstrating that pivotal sampling is more efficient than multinomial sampling. Indeed, it will follow from Proposition \ref{char:1:ops} in Section \ref{sec:ops} that it is sufficient to prove the result when sampling in $U_c$.

\begin{figure}[htb!]
    \setlength{\unitlength}{0.03cm}
    \begin{center}
    \begin{picture}(300,220)
    \thicklines
    % Le segment principal + Les bornes entières
    \put(0,180){\line(1,0){300}}
    \put(0,170){\line(0,1){20}}   \put(0,195){\makebox(0,0)[b]{$0$}}   \put(100,170){\line(0,1){20}} \put(100,195){\makebox(0,0)[b]{$1$}}
    \put(200,170){\line(0,1){20}} \put(200,195){\makebox(0,0)[b]{$2$}} \put(300,170){\line(0,1){20}} \put(300,195){\makebox(0,0)[b]{$3$}}
    % Les unités
    \put(0,160){\vector(1,0){28}}   \put(28,160){\vector(-1,0){28}} \put(14,145){\makebox(0,0)[b]{$\pi_{1}$}}
    \put(28,160){\vector(1,0){30}}   \put(58,160){\vector(-1,0){30}} \put(43,145){\makebox(0,0)[b]{$\pi_{2}$}}
    \put(60,160){\vector(1,0){68}}  \put(128,160){\vector(-1,0){68}} \put(95,145){\makebox(0,0)[b]{$\pi_{3}$}}
    \put(130,160){\vector(1,0){28}} \put(158,160){\vector(-1,0){28}} \put(145,145){\makebox(0,0)[b]{$\pi_{4}$}}
    \put(160,160){\vector(1,0){68}} \put(228,160){\vector(-1,0){68}} \put(195,145){\makebox(0,0)[b]{$\pi_{5}$}}
    \put(230,160){\vector(1,0){20}} \put(250,160){\vector(-1,0){20}} \put(238,145){\makebox(0,0)[b]{$\pi_{6}$}}
    \put(250,160){\vector(1,0){30}} \put(280,160){\vector(-1,0){30}} \put(265,145){\makebox(0,0)[b]{$\pi_{7}$}}
    \put(280,160){\vector(1,0){20}} \put(300,160){\vector(-1,0){20}} \put(288,145){\makebox(0,0)[b]{$\pi_{8}$}}

    \thinlines
    % Les a_i et b_i
    \put(60,185){\vector(1,0){39}}  \put(99,185){\vector(-1,0){39}} \put(79,187){\makebox(0,0)[b]{\tiny{$a_{1}$}}}
    \put(101,185){\vector(1,0){29}}  \put(130,185){\vector(-1,0){29}} \put(115,187){\makebox(0,0)[b]{\tiny{$b_{1}$}}}
    \put(160,185){\vector(1,0){39}}  \put(199,185){\vector(-1,0){39}} \put(179,187){\makebox(0,0)[b]{\tiny{$a_{2}$}}}
    \put(201,185){\vector(1,0){29}}  \put(230,185){\vector(-1,0){29}} \put(215,187){\makebox(0,0)[b]{\tiny{$b_{2}$}}}

    \thicklines
    % Le segment principal + Les bornes entières
    \put(0,80){\line(1,0){300}}
    \put(0,70){\line(0,1){20}}   \put(0,95){\makebox(0,0)[b]{$0$}}   \put(100,70){\line(0,1){20}} \put(100,95){\makebox(0,0)[b]{$1$}}
    \put(200,70){\line(0,1){20}} \put(200,95){\makebox(0,0)[b]{$2$}} \put(300,70){\line(0,1){20}} \put(300,95){\makebox(0,0)[b]{$3$}}
    % Les unités
    \put(0,60){\vector(1,0){58}}   \put(58,60){\vector(-1,0){58}} \put(30,45){\makebox(0,0)[b]{$\phi_{1}$}}
    \put(60,60){\vector(1,0){68}}  \put(128,60){\vector(-1,0){68}} \put(95,45){\makebox(0,0)[b]{$\phi_{2}$}}
    \put(130,60){\vector(1,0){28}} \put(158,60){\vector(-1,0){28}} \put(145,45){\makebox(0,0)[b]{$\phi_{3}$}}
    \put(160,60){\vector(1,0){68}} \put(228,60){\vector(-1,0){68}} \put(195,45){\makebox(0,0)[b]{$\phi_{4}$}}
    \put(230,60){\vector(1,0){70}} \put(300,60){\vector(-1,0){70}} \put(265,45){\makebox(0,0)[b]{$\phi_{5}$}}
    \thinlines
    % Les a_i et b_i
    \put(60,85){\vector(1,0){39}}  \put(99,85){\vector(-1,0){39}} \put(79,87){\makebox(0,0)[b]{\tiny{$a_{1}$}}}
    \put(101,85){\vector(1,0){29}}  \put(130,85){\vector(-1,0){29}} \put(115,87){\makebox(0,0)[b]{\tiny{$b_{1}$}}}
    \put(160,85){\vector(1,0){39}}  \put(199,85){\vector(-1,0){39}} \put(179,87){\makebox(0,0)[b]{\tiny{$a_{2}$}}}
    \put(201,85){\vector(1,0){29}}  \put(230,85){\vector(-1,0){29}} \put(215,87){\makebox(0,0)[b]{\tiny{$b_{2}$}}}
    \end{picture}
    \end{center}
    \vspace*{-2.2cm}\caption{A non-clustered population $U$ of size $8$ (top case) and the associated clustered population $U_c$ of size $5$ (bottom case) for a sample size $n=3$} \label{ops:Wi:Ji}
    \end{figure}
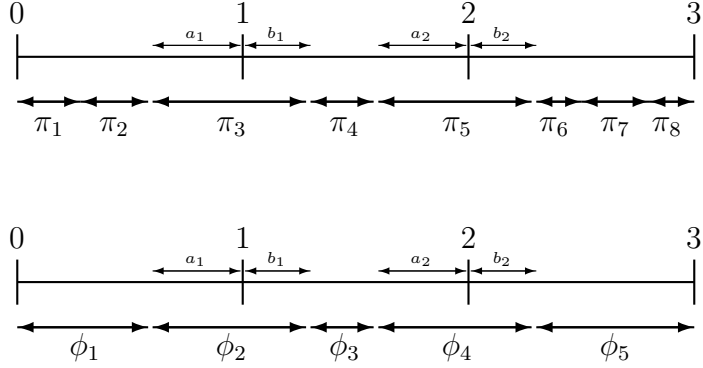

\noindent The population $U_c=\{u_1,\ldots,u_{2n-1}\}$ is obtained by clustering the units in $U$ as follows: each cross-border unit $k_i$ forms a separate cluster $u_{2i}$ of size $1$ with associated probability $\phi_{2i}=a_i+b_i$, for $i=1,\ldots,n-1$; the non cross-border units $k$ such that $k_{i-1}<k<k_i$ are grouped to form the cluster $u_{2i-1}$ with associated probability $\phi_{2i-1}=1-b_{i-1}-a_{i}$ for $i=1,\ldots,n$, where $k_0=0$ and $k_n=N+1$. We note $\phi=(\phi_1,\ldots,\phi_{2n-1})$. In the example presented in Figure \ref{ops:Wi:Ji}, the population $U$ contains two cross-border units $k_1=3$ and $k_2=5$. In the associated clustered population $U_c=\{u_1,\ldots,u_5\}$, the clusters of non cross-border units are $u_1$, which gathers units $k=1,2$; $u_3$, which contains the sole unit $k=4$; $u_5$, which gathers units $k=6,7,8$. The clusters of cross-border units are $u_2$ which contains the sole unit $k_1=3$, and $u_4$ which contains the sole unit $k_2=5$.

\section{Sampling Algorithms} \label{sec:ops}

\noindent Ordered pivotal sampling (ops) is recursively defined in Algorithm \ref{ops:pi:U}. For example, consider the population $U_c$ in Figure \ref{ops:Wi:Ji} with $\phi=(0.5,0.8,0.4,0.7,0.6)$, so that $a_1=0.5~$, $b_1=0.3$, $a_2=0.3$, $b_2=0.4$. Since $u_1$ is the sole non cross-border unit in $[0,1]$, we have $k_1=2$ and therefore, $H_1=u_1$. We take $(F_{1},L_{1})=(u_1,u_2)$ with probability $1-(1-b_1)^{-1} a_1=2/7$ and $(F_{1},L_{1})=(u_2,u_1)$ with probability $5/7$. In the first case, $H_2$ is selected among $\{u_2,u_3\}$ with probabilities proportional to $(0.3,0.4)$. If $H_2=u_3$, say, we take $(F_{2},L_{2})=(u_3,u_4)$ with probability $1-(1-b_2)^{-1} a_2=1/2$ and $(F_{2},L_{2})=(u_4,u_3)$ with probability $1/2$. In the second case, $H_3$ is selected among $\{u_3,u_5\}$ with probabilities proportional to $(0.4,0.6)$. This is also the last selected unit $F_3$. If $H_3=u_3$, say, the sample is $(F_1,F_2,F_3)=(u_1,u_4,u_3)$. The complete probability tree is given in the Supplementary Material.

\begin{algorithm}[htb!]
\begin{itemize}
\item One unit denoted as $H_1$ is selected among $\{1,\ldots,k_1-1\}$ with probabilities proportional to $(\pi_{1},\ldots,\pi_{k_1-1})$.
\item The unit $H_1$ faces the cross-border unit $k_1$. One unit, denoted as $F_1$, is selected in the sample while the other unit, denoted as $L_1$, goes on with the residual probability $b_1$. We have:
            \begin{eqnarray}
            (F_{1},L_{1}) & \equiv & \left\{ \begin{array}{ll}
                                 (H_{1},k_{1}) & \textrm{ with probability } 1-(1-b_1)^{-1} a_1,  \\
                                 (k_{1},H_{1}) & \textrm{ with probability } (1-b_1)^{-1} a_1.
                                 \end{array} \right.
            \end{eqnarray}
\item The $n-1$ remaining sampled units $\{F_{2},\ldots,F_n\}$ are drawn according to ordered pivotal sampling in the population $U^{(2)}=\{L_1,k_{1}+1,\ldots,N\}$ with inclusion probabilities $\pi^{(2)}=(b_1,\pi_{k_1+1},\ldots,\pi_{N})^{\top}$.
\item The final sample is $\{F_1,\ldots,F_n\}$.
\end{itemize}
\caption{Ordered pivotal sampling with parameter $\pi$ in $U$} \label{ops:pi:U}
\end{algorithm}

\noindent Multinomial sampling (ms) with parameter $\pi$ is a with replacement sampling algorithm, which consists of $n$ independent draws from the population $U$. At each draw, some unit $k$ in $U$ is selected with probability $n^{-1}\pi_k$. The variance of the Hansen-Hurvitz estimator under multinomial sampling is
    \begin{eqnarray} \label{var:han:hur:est:mul}
    V\left(\hat{t}_{y}^{ms}\right) & = & \sum_{k \in U} \pi_k \left(\check{y}_{k}-\frac{t_y}{n}\right)^2.
    \end{eqnarray}

\noindent Proposition \ref{char:1:ops} is a characterization of both sampling designs as two-stage procedures. The result for ordered pivotal sampling is given in \citet{cha:12}. The proof for multinomial sampling is omitted.

\begin{prop} \label{char:1:ops}
  Ordered pivotal sampling (respectively, multinomial sampling) with parameter $\pi$ in $U$ may be performed by two-stage sampling, with:
   \begin{enumerate}
     \item a first-stage selection of a sample $S_c^{ops}$ (respectively, $S_c^{ms}$) of $n$ clusters by means of ordered pivotal sampling (respectively, multinomial sampling) with parameter $\phi$ in the population $U_c$,
     \item an independent second-stage selection inside each $u_i \in S_c^{ops}$ (respectively, inside each $u_i \in S_c^{ms}$) of a sample $S_i$ of size $1$, with unit $k \in u_i$ selected with a probability $\phi_i^{-1} \pi_k$.
   \end{enumerate}
\end{prop}

\noindent From Proposition \ref{char:1:ops}, we have $\hat{t}_{y}^{ops} = \sum_{u_i \in S_c^{ops}} \sum_{k \in S_i} \check{y}_k$, which leads to
    \begin{eqnarray}
      V(\hat{t}_{y}^{ops}) & = & V \left\{ E (\hat{t}_{y}^{ops}|S_c^{ops}) \right\} + E \left\{ V (\hat{t}_{y}^{ops}|S_c^{ops}) \right\} \nonumber \\
                           & = & V \left\{ \sum_{u_i \in S_c^{ops}} \phi_i^{-1} \sum_{k \in u_i} y_k \right\} + E \left\{ \sum_{u_i \in S_c^{ops}} \sum_{k \in u_i} \phi_i^{-1} \pi_k \left( \check{y}_k -\check{Y}_i \right)^2 \right\} \nonumber \\
                           & = & V \left\{ \sum_{u_i \in S_c^{ops}} \check{Y}_i \right\} + \sum_{u_i \in U_c} \sum_{k \in u_i} \pi_k \left( \check{y}_k-\check{Y}_i \right)^2, \nonumber \\
                           & = & V \left( \hat{t}_{Y}^{ops}\right) + \sum_{u_i \in U_c} \sum_{k \in u_i} \pi_k \left(\check{y}_k-\check{Y}_i \right)^2, \label{var:2deg:ops}
    \end{eqnarray}
with $\hat{t}_{Y}^{ops}=\sum_{u_i \in S_c^{ops}} \check{Y}_i$, with $\check{Y}_i=\phi_i^{-1} Y_i$ and where $Y_i=\sum_{k \in u_i} y_k$. Similarly, we obtain from Proposition \ref{char:1:ops} that $\hat{t}_{y}^{ms} = \sum_{u_i \in S_c^{ms}} \sum_{k \in S_i} \check{y}_k$, which leads to
    \begin{eqnarray}
      V(\hat{t}_{y}^{ms}) = V \left( \hat{t}_{Y}^{ms}\right) + \sum_{u_i \in U_c} \sum_{k \in u_i} \pi_k \left(\check{y}_k-\check{Y}_i \right)^2
      & \textrm{ where } & \hat{t}_{Y}^{ms}=\sum_{u_i \in S_c^{ms}} \check{Y}_i.\label{var:2deg:ms}
    \end{eqnarray}
%    \begin{eqnarray}
%      V(\hat{t}_{y}^{ms}) & = & V \left( \hat{t}_{Y}^{ms}\right) + \sum_{u_i \in U_c} \sum_{k \in u_i} \pi_k \left(\check{y}_k-\check{Y}_i \right)^2 \label{var:2deg:ms}
%    \end{eqnarray}
%where $\hat{t}_{Y}^{ms}=\sum_{u_i \in S_c^{ms}} \check{Y}_i$.

\section{Comparison of the Sampling Algorithms} \label{sec:comp}

\noindent In equations (\ref{var:2deg:ops}) and (\ref{var:2deg:ms}), the first term on the right-hand side represents the variance due to the first stage of sampling, while the second term represents the variance due to the second stage of sampling. Clearly, ordered pivotal sampling and multinomial sampling may only differ with respect to the first term of variance. It is thus sufficient to prove that ordered pivotal sampling is more efficient when sampling in the clustered population $U_c$. The proof of Proposition \ref{prop:1} is available from the authors.

\begin{prop} \label{prop:1}
  We have:
    \begin{eqnarray}
    V\left(\hat{t}_{Y}^{ops}\right) & = & a_{1}(1-a_{1}-b_{1}) \left(\check{Y}_{1}-\check{Y}_{2}\right)^2 + E\left\{V\left(\left.\sum_{i=2}^n \check{Y}_{F_i} \right| F_{1} \right) \right\}, \label{var:recurs:hor:tho:est} \\
    V\left(\hat{t}_{Y}^{ms}\right) & \geq & a_{1}(1-a_{1}-b_{1}) \left(\check{Y}_{1}-\check{Y}_{2}\right)^2 + E\left\{V\left(\left.\sum_{u_i \in S_{(2)}^{ms}} \check{Y}_{i} \right|F_1 \right)\right\}, \label{var:recurs:han:hur:est}
    \end{eqnarray}
  with $S_{(2)}^{ms}$ a multinomial sample selected in $U_c^{(2)}=\{L_1,u_{3},\ldots,u_{2n-1}\}$ with parameter $\phi^{(2)}=(b_1,\phi_{3},\ldots,\phi_{2n-1})^{\top}$.
\end{prop}

\begin{theo} \label{sup:ops}
  Ordered pivotal sampling with parameter $\pi$ is more accurate than Multinomial sampling with parameter $\pi$.
\end{theo}

\noindent Theorem \ref{sup:ops} follows from Proposition \ref{prop:1} with a proof by induction. It implies the weak consistency of the HT-estimator, as summarized in Theorem \ref{mean:squa:cons}.
\begin{theo} \label{mean:squa:cons}
  Assume that: \\
  H1: there exists some constant $A_1$ s.t. $\sum_{k \in U} \pi_k \left(\check{y}_k-n^{-1} t_y \right)^2 \leq A_1 N^2 n^{-1}$. \\
  Then $E\left\{ N^{-1} (\hat{t}_{y}^{ops}-t_y) \right\}^2 = O(n^{-1})$, and the HT-estimator is weakly consistent for $t_y$.
\end{theo}
\noindent The proof is straightforward. The moment assumption H1 will hold in particular: if there exists some constant $A_2>0$ such that $\min_{k \in U} \pi_k \geq A_2 N^{-1} n$; and if there exists some constant $A_3$ such that $N^{-1} \sum_{k \in U} y_k^2 \leq A_3$.

\noindent A drawback of ordered pivotal sampling lies in variance estimation. The second-order inclusion probabilities can be computed exactly \citep{cha:12}, but many of them are usually equal to $0$ which results in a biased variance estimator. Denote by $v_{HH}(\hat{t}_y^{ops}) = (n-1)^{-1} n \sum_{k \in S^{ops}} \left( \check{y}_k - n^{-1} \hat{t}_y^{ops}\right)^2$ the HH-variance estimator applied to pivotal sampling. From $v_{HH}(\hat{t}_y^{ops}) = n(n-1)^{-1} \sum_{k \in S^{ops}} \left( \check{y}_k - \frac{t_y}{n}\right)^2 - (n-1)^{-1} \left(\hat{t}_y^{ops}-t_y\right)^2$, we have
    \begin{eqnarray} \label{theo4:eq4}
      E\left[v_{HH}(\hat{t}_y^{ops})-V(\hat{t}_y^{ops}) \right] & = & \frac{n}{n-1} \left[V(\hat{t}_y^{ms})-V(\hat{t}_y^{ops}) \right].
    \end{eqnarray}
It follows from (\ref{theo4:eq4}) and Theorem \ref{sup:ops} that the HH-variance estimator can always be used as a conservative variance estimator for pivotal sampling. This result is particularly of interest in a spatial sampling context \citep{gra:lun:sch:12} when the joint selection of neighbouring units is avoided so as to build an efficient sampling design.

\noindent Theorem \ref{sup:ops} can be easily extended to any randomized version of ordered pivotal sampling. Denote by $\sigma$ a random permutation of the units in $U$. Randomized pivotal sampling is obtained by applying Algorithm \ref{ops:pi:U} to the randomized population $U_{\sigma}=\{\sigma(1),\ldots,\sigma(N)\}$ with parameter ${\pi}_{\sigma}=(\pi_{\sigma(1)},\ldots,\pi_{\sigma(N)})^{\top}$. Then it is easily shown that randomized pivotal sampling with parameter $\pi$ is more accurate than multinomial sampling with parameter $\pi$. This implies that under randomized pivotal sampling, the HT-estimator is weakly consistent for $t_y$ under the assumption H1 and the HH-variance estimator is always conservative for the true variance.

\section{Simulation study} \label{sec:simu}

\noindent We conducted a simulation study to confirm our theoretical results. We used the clustered population $U_c$ associated to the sample size $n=3$, presented in Figure \ref{ops:Wi:Ji}. We considered all the possible sets of inclusion probabilities with a skip of $0.05$; that is, all the possible sets of inclusion probabilities $\phi=(\phi_1,\phi_2,\phi_3,\phi_4,\phi_5)$ such that for any $i=1,\ldots,5$, $\phi_i=0.05 m_i$ for some integer $m_i$, with $0<\phi_i<1$ and $\sum_{i=1}^5 \phi_i = 3$. This led to $24,396$ cases.

\noindent As proved in \citet{gab:84}, a sampling design is more efficient than multinomial sampling iff the matrix $B=(\phi_{ij}/\phi_j)_{i,j \in U_c}$ has its second largest eigenvalue $\lambda_2 \leq 1$, with $\phi_{ij}$ the second-order inclusion probability for clusters $u_i$ and $u_j$. In fact, $\lambda_2$ corresponds to the largest possible value for the ratio of the variances under ordered pivotal sampling and under multinomial sampling, see \citet[page 69, equation (3)]{gab:90}.

\noindent For each of the $24,396$ cases, we computed $B$ using the formulas for second-order inclusion probabilities given in Theorem 5.2 in \citet{cha:12}, and the second largest eigenvalue of $B$. The values of $\lambda_2$ ranged from $0.625$ to $0.991$, confirming the result. We conducted a similar simulation study on the clustered population $U_c$ for a sample size $n=5$, with a skip of $0.10$. For each of the $31,998$ cases, the values of $\lambda_2$ ranged from $0.666$ to $0.975$.

\end{document}